\documentclass[12pt]{article}
\def\Z{\mathbb Z}

\def\N{\mathbb N}
\def\A{\mathcal A}

\def\C{\mathcal C}

\def\pf{\begin{proof}}
\def\pfk{\end{proof}}

\newcommand{\coloneq}{\mathrel{\mathop:}=}

\usepackage{amssymb,amsmath,amsthm,amsfonts,amscd}
\usepackage[T1]{fontenc}

\newtheorem{lem}{Lemma}[section]
\newtheorem{prop}[lem]{Proposition}
\newtheorem{coro}[lem]{Corollary}
\newtheorem{thm}[lem]{Theorem}
\theoremstyle{definition}
\newtheorem{de}[lem]{Definition}
\newtheorem{pozn}[lem]{Remark}
\newtheorem{ex}[lem]{Example}

\begin{document}
\title{Palindromic complexity of infinite words associated with simple Parry numbers}

\author{Petr Ambro\v z$^{(1)}$$^{(2)}$ \and
Christiane Frougny$^{(2)}$$^{(3)}$
\and Zuzana Mas\'akov\'a$^{(1)}$
\and Edita Pelantov\'a$^{(1)}$}

\maketitle
\thispagestyle{empty}

\begin{center}
(1) Doppler Institute for Mathematical Physics and Applied Mathematics and\\
Department of Mathematics, FNSPE, Czech Technical University, \\
Trojanova 13, 120 00
Praha 2, Czech Republic\\
{\tt Ampy@linux.fjfi.cvut.cz, Masakova@km1.fjfi.cvut.cz, Pelantova@km1.fjfi.cvut.cz}\\
(2) LIAFA, UMR 7089 CNRS \& Universit\'e Paris 7\\
2 place Jussieu, 75251 Paris Cedex 05, France \\
(3) Universit\'e Paris 8\\
{\tt Christiane.Frougny@liafa.jussieu.fr}
\end{center}

\begin{abstract}
A simple Parry number is a real number $\beta>1$ such that the R\'enyi
expansion of $1$ is finite, of the form $d_\beta(1)=t_1 \cdots t_m$.
We study the palindromic structure of infinite aperiodic words $u_\beta$ that are the fixed point
of a substitution associated with a simple Parry number $\beta$.
It is shown that the word $u_\beta$ contains infinitely many palindromes
if
and only if $t_1=t_2= \cdots=t_{m-1}\ge t_m$.
Numbers $\beta$ satisfying this condition are the so-called {\em confluent} Pisot numbers.
If $t_m=1$ then $u_\beta$ is an Arnoux-Rauzy word.
We show that if $\beta$ is a confluent Pisot number
then $
{\mathcal P}(n+1)+ {\mathcal P}(n) = {\mathcal C}(n+1) - {\mathcal C}(n)
+2$, where ${\mathcal P}(n)$ is the number of palindromes
and ${\mathcal C}(n)$ is the number of factors
of length $n$ in $u_\beta$.
We then
give a complete description of the set of palindromes, its structure and properties.
\end{abstract}

\section{Introduction}

Infinite aperiodic words over a finite alphabet are suitable models for one-dimensional
quasicrystals, i.e.\ non-crystallographic materials displaying long-range order, since
they define one-dimensional Delaunay sets with finite local complexity.
The first quasicrystal was discovered
in 1984: it is a solid structure presenting a local symmetry of
order 5, i.e.\ a local invariance under rotation of $\pi/5$,
and it is linked to the golden ratio and to the Fibonacci
substitution. The Fibonacci substitution, given by
$$ 0 \mapsto 01, \; 1 \mapsto 0 ,$$
defines a quasiperiodic selfsimilar tiling of the positive real line, and
is a historical model of a one-dimensional quasicrystal.
The fixed point of the substitution is the infinite word
$$010010101\cdots$$
The description and the properties of this tiling
use a number system in base the golden ratio.

A more general theory has been elaborated with Pisot
numbers\footnote{A {\em Pisot number} is an algebraic integer
$>1$ such that the other roots of its minimal polynomial have a
modulus less than $1$. The golden ratio and the natural integers
are Pisot numbers.} for base, see \cite{Ga,bfgk98}. Note that so
far, all the quasicrystals discovered by physicists present local
symmetry of order 5 or 10, 8, and 12, and are modelized using
quadratic Pisot units, namely the golden ratio  for order 5 or 10, $1+\sqrt{2}$
for order 8, and $2+\sqrt{3}$ for order 12.

For the description of physical properties of these materials it is important to know
the combinatorial properties of the infinite aperiodic words, such as the factor
complexity, which corresponds to the number of local configurations of atoms in the
material, or the palindromic structure of the aperiodic words, describing local symmetry
of the material. The palindromic structure of the infinite words has been proven
important for the description of the spectra of Schr\"odinger operators with potentials
adapted to aperiodic structures~\cite{HKS}.

The most studied infinite aperiodic word is the Fibonacci word, which is the paradigm
of the notion of
Sturmian words. Sturmian words are binary aperiodic words with minimal factor
complexity, i.e.\ ${\mathcal C}(n)=n+1$ for $n\in\N = \{0, 1, 2, \ldots \}$. There exist several equivalent
definitions of Sturmian words see~\cite{berstel}, or~\cite[Chapter 2]{lothaire2}. From our point of view the
characterization of Sturmian words using palindromes~\cite{DrPi} is particularly
interesting.

Sturmian words can be generalized in several different ways to words over $m$-letter
alphabet, namely to Arnoux-Rauzy words of order $m$, see \cite{AR,berstel},
or to infinite words coding
$m$-interval exchange \cite{keane,rauzy}. The Sturmian case is included for $m=2$.

Arnoux-Rauzy words and words coding generic $m$-interval exchange have factor complexity
${\mathcal C}(n)=(m-1)n+1$ for $n\in\N$, \cite{keane}. For Arnoux-Rauzy words the palindromic structure
is also known~\cite{JuPi,Da}: for every $n$ the number ${\mathcal P}(n)$
of palindromes of length $n$
is equal to ${\mathcal P}(n)=1$ if $n$ is even and to ${\mathcal P}(n)=m$
if $n$ is odd. The
palindromic structure of infinite words coding $m$-interval exchange is more complicated.
The existence of palindromes of arbitrary length depends on the permutation which
exchanges the intervals. For $m=3$ and the permutation $\pi=(321)$ the result is given
in~\cite{Da}, for general $m$ in~\cite{BaPe}.

As we have seen for the Fibonacci word, infinite aperiodic words can also be obtained as the fixed point of a
 substitution canonically associated with a number system
where the base is an irrational number $\beta$, the so-called
$\beta$-expansions
introduced by R\'enyi~\cite{renyi}. The words $u_\beta$ are defined in the case that
 $\beta$ is a Parry number, that is to say when the R\'enyi expansion of
1 is eventually periodic or finite,
 see Section~\ref{prel} for definitions. These words provide a good model
of one-dimensional
quasicrystals~\cite{Ga}.
The factor
complexity of these words is at most linear, because they are fixed points of primitive
substitutions~\cite{quef}. The exact values of the complexity function ${\mathcal C}(n)$
for a large class of Parry numbers $\beta$ can be found in~\cite{FrMaPe} and some
partial results about other Parry numbers $\beta$ are to be found in~\cite{FrMaPe2}.

\medskip

This paper is devoted to the description of the palindromic structure of the infinite
words $u_\beta$, when $\beta$ is a simple Parry number, with the R\'enyi expansion of
1 being of the form $d_\beta(1)=t_1 \cdots t_m$.
We first show that the word $u_\beta$ contains infinitely many palindromes
if
and only if $t_1=t_2= \cdots=t_{m-1}\ge t_m$.
Numbers $\beta$ satisfying this condition have been introduced and studied in~\cite{frougny} from the point of view
of linear numeration systems. Confluent linear numeration systems are exactly
those for which there is no
propagation of the carry to the right in the process of normalization, which
consists of transforming a non-admissible representation on the canonical alphabet
of a number into the
admissible $\beta$-expansion of that number.
Such a number $\beta$ is known to be a Pisot number,
and will be called a {\em confluent} Pisot number.
We also know from \cite{FrMaPe} that the infinite word $u_\beta$ is an Arnoux-Rauzy sequence
if and only if it is a confluent Pisot number with the last
coefficient $t_m$ being equal to $1$; then $\beta$ is an algebraic unit.

In the sequel $\beta$ is a confluent Pisot number.
We then determine the palindromic complexity, that
is ${\mathcal P}(n)$, the
number of palindromes in $u_\beta$ of length $n$.
In the description of ${\mathcal P}(n)$
we use the notions introduced in~\cite{FrMaPe} for the factor complexity. The connection
of the factor and palindromic complexity is not surprising. For example, in~\cite{ABCD}
the authors give an upper estimate of the palindromic complexity ${\mathcal P}(n)$ in
terms of ${\mathcal C}(n)$.

In this paper we show that if the length of palindromes is not bounded, which is equivalent to
$\limsup_{n\to\infty}{\mathcal P}(n)>0$, then
\begin{equation}\label{eq:PaC}
{\mathcal P}(n+1)+ {\mathcal P}(n) = {\mathcal C}(n+1) - {\mathcal C}(n)
+2\,,\qquad\hbox{ for }\ n\in\N\,.
\end{equation}

In general it is has been shown~\cite{BaMaPe2} that for a uniformly recurrent word with
$\limsup_{n\rightarrow\infty}{\mathcal P}(u)>0$ the inequality
$$
{\mathcal P}(n+1)+ {\mathcal P}(n) \leq {\mathcal C}(n+1) - {\mathcal C}(n)
+2
$$
holds for all $n\in\N$. Moreover, the authors proved the
formula~\eqref{eq:PaC} to be valid for infinite words coding the
$r$-interval exchange. Finally, it is known that the
formula~\eqref{eq:PaC} holds also for Arnoux-Rauzy
words~\cite{ABCD} and for complementation-symmetric
sequences~\cite{Da}.

We then give a complete description of the set of palindromes, its structure and properties.
The exact palindromic complexity of the word $u_\beta$ is given in Theorem~\ref{thm:10}.

Further on, we study the occurrence of palindromes of arbitrary length in the prefixes of
the word $u_\beta$. It is known~\cite{DrJuPi} that every word $w$ of length $n$ contains at most $n+1$
different palindromes. The value by which the number of palindromes differs from $n+1$
is called the {\em defect} of the word $w$. Infinite words whose every prefix has defect 0 are
called {\em full}. We show that whenever $\limsup_{n\to\infty}{\mathcal P}(n)>0$, the infinite
word $u_\beta$ is full.

\section{Preliminaries}\label{prel}

Let us first recall the basic notions which we work with, for more
details reader is referred to~\cite{lothaire2}. An {\em alphabet}
is a finite set whose elements are called letters. A {\em finite
word} $w=w_1w_2\cdots w_n$ on the alphabet ${\mathcal A}$ is a
concatenation of letters. The length $n$ of the word $w$ is
denoted by $|w|$. The set of all finite words together with the
empty word $\varepsilon$ equipped with the operation of
concatenation is a free monoid over the alphabet $\A$, denoted by
$\A^*$.

An infinite sequence of letters of $\A$ of the form
$$
u_0u_1u_2\cdots,\qquad \cdots u_2u_1u_0,\qquad\hbox{ or }\qquad \cdots
u_{-2}u_{-1}u_0u_1u_2 \cdots
$$
is called {\em right infinite} word, {\em left infinite} word, or
{\em two-sided infinite} word, respectively. If for a two-sided
infinite word the position of the letter indexed by 0 is
important, we introduce {\em pointed} two-sided infinite words,
$\cdots u_{-2}u_{-1}|u_0u_1u_2\cdots$.

A {\em factor} of a word $v$ (finite or infinite) is a finite word
$w$ such that there exist words $v_1,v_2$ satisfying $v=v_1wv_2$.
If $v_1=\varepsilon$, then $w$ is called a {\em prefix} of $v$, if
$v_2=\varepsilon$, then $w$ is a {\em suffix} of $v$. For a finite
word $w=w_1w_2\cdots w_n$ with a prefix $v=w_1\cdots w_k$, $k\leq
n$, we define $v^{-1}w\coloneq w_{k+1}\cdots w_n$.

On the set $\A^*$ we can define the operation $\sim$ which to a finite word $w=w_1\cdots
w_n$ associates $\widetilde{w}=w_n\cdots w_1$. The word $\widetilde{w}$ is called the
{\em reversal} of $w$. A finite word $w\in\A^*$ for which $w=\widetilde{w}$ is called a
{\em palindrome}.

The set of all factors of an infinite word $u$ is called the {\em language} of $u$ and denoted
by ${\mathcal L}(u)$. The set of all palindromes in ${\mathcal L}(u)$ is denoted by
${\mathcal Pal}(u)$. The set of words of length $n$ in ${\mathcal L}(u)$, respectively in
${\mathcal Pal}(u)$ determines the {\em factor}, respectively {\em palindromic complexity} of the
infinite word $u$. Formally, the functions ${\mathcal C}:\N\to\N$, ${\mathcal
P}:\N\to\N$ are defined by
\begin{eqnarray*}
{\mathcal C}(n)&\coloneq& \#\,\{w \mid w\in{\mathcal L}(u),\, |w|=n\}\,,\\
{\mathcal P}(n)&\coloneq& \#\,\{w \mid w\in{\mathcal Pal}(u),\, |w|=n\}\,.
\end{eqnarray*}
Obviously, we have ${\mathcal P}(n)\leq {\mathcal C}(n)$ for $n\in\N$. We have moreover
$$
{\mathcal P}(n)\leq \frac{16}{n} \,
{\mathcal C}\left(n+\Big\lfloor\frac{n}{4}\Big\rfloor\right)\,,
$$
as shown in~\cite{ABCD}.

For the determination of the factor complexity important is the notion of the so-called
left or right special factors, introduced in~\cite{cassaigne}.
The extension of a factor $w\in{\mathcal L}(u)$ by a
letter to the left is called the {\em left extension} of $w$, analogously we define the {\em right
extension} of a factor $w$. Formally, we have the sets
\begin{eqnarray*}
{\rm Lext}(w)&\coloneq& \{aw \mid aw\in{\mathcal L}(u)\}\,,\\
{\rm Rext}(w)&\coloneq& \{wa \mid wa\in{\mathcal L}(u)\}\,.
\end{eqnarray*}
If $\#{\rm Lext}(w)\geq 2$, we say that $w$ is a {\em left special factor} of the infinite word
$u$. Similarly, if $\#{\rm Rext}(w)\geq 2$, then $w$ is a {\em right special factor}
 of $u$.
For the first difference of complexity we have
$$
\Delta {\mathcal C}(n) \ = \ {\mathcal C}(n+1) - {\mathcal C}(n) \ = \
\sum_{w\in{\mathcal L}(u), \ |w|=n } (\#{\rm Lext}(w)-1)\,.
$$
In this formula we can exchange ${\rm Lext}(w)$ with ${\rm Rext}(w)$.

Infinite words which have for each $n$ at most one left special
factor and at most one right special factor are called {\em
episturmian} words~\cite{JuPi}. Arnoux-Rauzy words of order $m$ are special
cases of episturmian words; they are defined as words on a
$m$-letter alphabet such that for every $n$ there exist exactly
one left special factor $w_1$ and exactly one right special factor
$w_2$. Moreover, these special factors satisfy $\#{\rm
Lext}(w_1)=\#{\rm Rext}(w_2)=m$.

Analogically to the case of factor complexity, for the palindromic complexity it is
important to define the palindromic extension: If for a palindrome $p\in{\mathcal
Pal}(u)$ there exists a letter $a$ such that $apa\in{\mathcal Pal}(u)$, then we call the
word $apa$ the {\em palindromic extension} of $p$.

A mapping on a free monoid $\A^*$ is called a {\em morphism} if
$\varphi(vw)=\varphi(v)\varphi(w)$ for all $v,w\in\A^*$. Obviously, for determining the
morphism it is sufficient to define $\varphi(a)$ for all $a\in\A$. The action of a morphism
can be naturally extended on right infinite words by the prescription
$$
\varphi(u_0u_1u_2\cdots)\coloneq \varphi(u_0)\varphi(u_1)\varphi(u_2)\cdots\,.
$$
A non-erasing\footnote{A morphism $\varphi$ on an alphabet $\A$ is \emph{non-erasing} if for
any $a\in\A$ the image $\varphi(a)$ is a non-empty word.}
morphism $\varphi$, for which there exists a letter $a\in\A$ such that
$\varphi(a)=aw$ for some non-empty word $w\in\A^*$, is called a {\em substitution}. An infinite
word $v$ such that $\varphi(v)=v$ is called a {\em fixed point} of the substitution $\varphi$.
Obviously, any substitution has at least one fixed point,
namely $\lim_{n\to\infty} \varphi^n(a)$. Assume that there exists an index $k\in\N$ such
that for every pair of letters $i,j\in\A$ the word $\varphi^k(i)$ contains as a factor
the letter $j$. Then the substitution $\varphi$ is called {\em primitive}.

Similarly, one can extend the action of a morphism to left infinite words.
For a pointed two-sided infinite word $u = \cdots u_{-3}u_{-2}u_{-1}|u_0u_1\cdots$
we define action of a morphism $\varphi$ by
$\varphi(u) = \cdots \varphi(u_{-3}) \varphi(u_{-2}) \varphi(u_{-1})| \varphi(u_0) \varphi(u_1)\cdots$.
One can also define analogically the notion of a fixed point.

Right infinite words which will be studied in this paper, are connected with the
R\'enyi $\beta$-expansion of real numbers \cite{renyi}. For a real number $\beta>1$ the
transformation $T_\beta:[0,1]\to[0,1)$ is defined by the prescription
$$
T_\beta(x) \coloneq \beta x- \lfloor\beta x\rfloor\,.
$$
The sequence of non-negative integers $(t_n)_{n\geq 1}$ defined by
$t_i=\lfloor\beta T^{i-1}(1)\rfloor$ satisfies
$1=\frac{t_1}{\beta} + \frac{t_2}{\beta^2} +
\frac{t_3}{\beta^3}+\cdots$. It is called the {\em R\'enyi
expansion} of $1$ and denoted by
$$
d_\beta(1)=t_1t_2t_3\cdots\,.
$$
In order that the  sequence $t_1t_2t_3\cdots$ be the R\'enyi expansion of 1 for some
$\beta$, it must satisfy the so-called Parry condition \cite{parry}
$$
t_it_{i+1}t_{i+2}\cdots \ \prec \ t_1t_2t_3\cdots \qquad\hbox{ for all }\ i=2,3,\ldots\,,
$$
where the symbol $\prec$ stands for "lexicographically strictly
smaller". A number $\beta>1$ for which $d_\beta(1)$ is eventually
periodic is called a {\em Parry number}. If moreover $d_\beta(1)$
has only finitely many non-zero elements, we say that $\beta$ is a
{\em simple} Parry number and in the notation for $d_\beta$ we
omit the ending zeros, i.e.\ $d_\beta(1)=t_1t_2\cdots t_m$, where
$t_m\neq 0$.

A {\em Pisot number} is an algebraic integer such that all its Galois conjugates are in modulus less
than $1$. A  Pisot number is a Parry number \cite{bertrand}.
It is known that a Parry number is a {\em Perron number}, i.e.\ an algebraic integer all of whose
conjugates are in modulus less than $\beta$. Solomyak~\cite{solo} has shown that all
conjugates of a Parry number lie inside the disc of radius $\frac12(1+\sqrt5)$, i.e.\ the
golden ratio.

With every Parry number one associates a canonical substitution $\varphi_\beta$,
see~\cite{fabre}. For a simple Parry number $\beta$ with $d_\beta(1)=t_1t_2\cdots t_m$
the substitution $\varphi=\varphi_\beta$ is defined on the alphabet
$\A=\{0,1,\ldots,m-1\}$ by
\begin{equation}\label{e:substi}
\begin{array}{ccl}
\varphi(0)&=&0^{t_1}1\\
\varphi(1)&=&0^{t_2}2\\
&\vdots&\\
\varphi(m-2)&=&0^{t_{m-1}}(m-1)\\
\varphi(m-1)&=&0^{t_m}
\end{array}
\end{equation}
The notation $0^k$ in the above stands for a concatenation of $k$ zeros. The substitution
$\varphi$ has a unique fixed point, namely the word
$$
u_\beta \coloneq \lim_{n\to\infty} \varphi^n(0)\,,
$$
which is the subject of the study of this paper. The
substitution~\eqref{e:substi} is primitive, and thus according
to~\cite{quef}, the factor complexity of its fixed point is
sublinear. The exact values of ${\mathcal C}(n)$
for $u_{\beta}$ with $d_{\beta}=t_1\cdots t_m$ satisfying
$t_1 > \max \{t_2,\ldots, t_{m-1}\}$ or $t_1 = t_2 = \cdots = t_{m-1}$
can be found in~\cite{FrMaPe}. The determination of the palindromic complexity
of $u_\beta$ is the aim of this article.

A similar canonical substitution is defined for non-simple Parry numbers.
Partial results about the factor and palindromic complexity of $u_{\beta}$
for non-simple Parry numbers $\beta$ can be found in~\cite{Lubka,FrMaPe2}.

One can define the canonical substitution $\varphi_{\beta}$ even if
the R\'enyi expansion $d_\beta(1)$ is infinite non-periodic, i.e.\ $\beta$ is not a
Parry number. In this case, however, the substitution and its fixed point are
defined over an infinite alphabet. The study of such words $u_{\beta}$ is
out of the scope of this paper.

\section{Words $u_\beta$ with bounded number of palindromes}

The infinite word $u_\beta$ associated with a Parry number $\beta$ is a fixed point of a
primitive substitution. This implies that the word $u_\beta$ is uniformly recurrent~\cite{fogg}.
Let us recall that an infinite word $u$ is called {\em uniformly recurrent} if every factor
$w$ in ${\mathcal L}(u)$ occurs in $u$ with bounded gaps.

\begin{lem}
If the language ${\mathcal L}(u)$ of a uniformly recurrent word $u$ contains infinitely
many palindromes, then ${\mathcal L}(u)$ is closed under reversal.
\end{lem}

\pf
 From the definition of a uniformly recurrent word $u$ it follows that for every $n\in\N$
 there exists an integer $R(n)$ such that every arbitrary factor of $u$ of length $R(n)$
 contains all factors of $u$ of length $n$. Since we assume that ${\mathcal Pal}(u)$ is
 an infinite set, it must contain a palindrome $p$ of length $\geq R(n)$. Since $p$
 contains all factors of $u$ of length $n$, and $p$ is a palindrome, it contains with every $w$ such that
 $|w|=n$ also its reversal $\widetilde{w}$. Thus $\widetilde{w}\in{\mathcal L}(u)$.
This consideration if valid for all $n$ and thus the statement of the lemma is proved.
 \pfk

Note that this result
was first stated, without proof, in~\cite{DrJuPi}.

The fact that the language is closed under reversal is thus a necessary condition so that
a uniformly recurrent word has infinitely many palindromes. The converse is not true \cite{bbcf}.

For infinite words $u_\beta$ associated with simple Parry numbers $\beta$ the invariance
of ${\mathcal L}(u_\beta)$ under reversal is studied in~\cite{FrMaPe}.

\begin{prop}[\cite{FrMaPe}]
Let $\beta>1$ be a simple Parry number such that $d_\beta(1)=t_1t_2\cdots t_m$\,.

\begin{enumerate}
\item
The language ${\mathcal L}(u_\beta)$ is closed under reversal, if and only if
$$
{\rm Condition\ (C):} \hspace*{2.5cm}t_1=t_2=\cdots = t_{m-1}\,. \hspace*{3cm}
$$

\item
The infinite word $u_\beta$ is an Arnoux-Rauzy word if and only if Condition (C) is
satisfied and $t_m=1$.
\end{enumerate}
\end{prop}
\begin{coro}
Let $\beta$ be a simple Parry number which does not satisfy Condition (C). Then there
exists $n_0\in\N$ such that  ${\mathcal P}(n)=0$ for $n\geq n_0$.
\end{coro}

Numbers $\beta$ satisfying Condition (C) have been introduced and studied in~\cite{frougny} from the point of view
of linear numeration systems. Confluent linear numeration systems are exactly
those for which there is no
propagation of the carry to the right in the process of normalization, which
consists of transforming a non-admissible representation on the canonical alphabet
of a number into the
admissible $\beta$-expansion of that number.
A number $\beta$ satisfying Condition (C) is known to be a Pisot number,
and will be called a {\em confluent} Pisot number.

Set
$$
t\coloneq t_1=t_2=\cdots = t_{m-1}\qquad\hbox{ and }\qquad s\coloneq t_m\,.
$$
{}From the Parry condition for the R\'enyi expansion of $1$ it follows that $t\geq s\geq
1$. Then the substitution $\varphi$ is of the form
\begin{equation}\label{e:1}
\begin{array}{ccl}
\varphi(0)&=&0^t1\\
\varphi(1)&=&0^t2\\
&\vdots&\\
\varphi(m-2)&=&0^t(m-1)\\
\varphi(m-1)&=&0^s
\end{array}
\qquad \qquad t\geq s\geq 1\,.
\end{equation}
Note that in the case $s=1$, the number $\beta$ is an algebraic
unit, and the corresponding word $u_\beta$ is an Arnoux-Rauzy
word, for which the palindromic complexity is known. Therefore in
the paper we often treat separately the cases $s\geq 2$ and $s=1$.

\section{Palindromic extensions in $u_\beta$}

In the remaining part of the paper we study the palindromic structure of the words
$u_\beta$ for confluent Pisot numbers $\beta$.

For an Arnoux-Rauzy word $u$ (and thus also for a Sturmian word)
it has been shown that for every palindrome $p\in{\mathcal L}(u)$
there is exactly one letter $a$ in the alphabet, such that
$apa\in{\mathcal L}(u)$, i.e.\ any palindrome in an Arnoux-Rauzy
word has exactly one palindromic extension~\cite{Da}. Since the
length of the palindromic extension $apa$ of $p$ is $|apa|=|p|+2$,
we have for Arnoux-Rauzy words ${\mathcal P}(n+2)={\mathcal P}(n)$
and therefore
$$
{\mathcal P}(2n)={\mathcal P}(0)=1\qquad\hbox{ and }\qquad {\mathcal P}(2n+1)={\mathcal
P}(1)=\#\A\,.
$$
Determining the number of palindromic extensions for a given palindrome of $u_\beta$ is
essential also for our considerations here. However, let us first introduce the following
notion.

\begin{de}
We say that a palindrome $p_1$ is a {\em central factor} of a palindrome $p_2$ if there exists a
finite word $w\in{\mathcal A}^*$ such that $p_2=wp_1\widetilde{w}$.
\end{de}

For example, a palindrome is a central factor of its palindromic extensions.

The following simple result can be easily obtained from the
form of the substitution~\eqref{e:1}, and is
a special case of a result given in~\cite{FrMaPe}.

\begin{lem}[\cite{FrMaPe}]\label{l:1}
All factors of $u_\beta$ of the form $X0^nY$ for $X,Y\neq 0$ are the following
\begin{equation}\label{e:2}
X0^t1, \ 10^tX\hbox{ with } X\in\{1,2,\ldots,m-1\},\ \hbox{ and }\ 10^{t+s}1\,.
\end{equation}
\end{lem}

\begin{pozn}\label{pozn:2.3}\
\begin{enumerate}
\item
Every pair of non-zero letters in $u_\beta$ is separated by a
block of at least $t$ zeros. Therefore every palindrome
$p\in{\mathcal L}(u_\beta)$ is a central factor of a palindrome
with prefix and suffix $0^t$.

\item
Since $\varphi({\mathcal A})$ is a suffix code, the coding given by the substitution
$\varphi$ is uniquely decodable. In particular, if $w_1\in{\mathcal L}(u_\beta)$ is a
factor with the first and the last letter non-zero, then there exist a factor
$w_2\in{\mathcal L}(u_\beta)$ such that $0^tw_1=\varphi(w_2)$.
\end{enumerate}
\end{pozn}

\begin{prop}\label{p:1}\
\begin{itemize}
\item[(i)] Let $p\in{\mathcal L}(u_\beta)$. Then $p\in{\mathcal Pal}(u_\beta)$ if and only if
$\varphi(p)0^t\in{\mathcal Pal}(u_\beta)$.

\item[(ii)] Let $p\in{\mathcal Pal}(u_\beta)$. The number of palindromic extensions of $p$ and
$\varphi(p)0^t$ is the same, i.e.\
$$ \#\{a\in{\mathcal A}\mid apa
\in {\mathcal Pal}(u_\beta)\} = \#\{a\in{\mathcal A}\mid
a\varphi(p)0^ta \in {\mathcal Pal}(u_\beta)\}\,.
$$
\end{itemize}
\end{prop}

\pf (i)  Let $p=w_0w_1\cdots w_{n-1}\in{\cal L}(u_\beta)$. Let us study under which
conditions the word $\varphi(p)0^t$ is also a palindrome, i.e.\ when
\begin{equation}\label{e:3}
\varphi(w_0)\varphi(w_1)\cdots\varphi(w_{n-1})0^t = 0^t\widetilde{\varphi(w_{n-1})}\cdots
\widetilde{\varphi(w_1)}\widetilde{\varphi(w_0)}\,.
\end{equation}
The substitution $\varphi$ has the property that for each letter $a\in{\cal A}$ it
satisfies $\widetilde{\varphi(a)}=0^{-t}\varphi(a)0^t$. Using this property, the
equality~\eqref{e:3} can be equivalently written as
$$
\varphi(p)=\varphi(w_0)\cdots\varphi(w_{n-1})=\varphi(w_{n-1})\cdots\varphi(w_0) =
\varphi(\widetilde{p})\,.
$$
As a consequence of unique decodability of $\varphi$ we obtain that~\eqref{e:3} is valid
if and only if $p=\widetilde{p}$.

(ii) We show that for a palindrome $p$ it holds that
$$
apa\in{\mathcal Pal}(u_\beta) \quad\iff\quad b\varphi(p)0^tb \in{\mathcal
Pal}(u_\beta)\,,\ \hbox{ where } b\equiv a+1 \pmod m\,,
$$
which already implies the equality of the number of palindromic extensions of palindromes
$p$ and $\varphi(p)0^t$.

Let $apa\in{\mathcal Pal}(u_\beta)$. Then
$$
\varphi(a)\varphi(p)\varphi(a)0^t = \left\{\begin{array}{ll}
 0^t(a+1)\varphi(p)0^t(a+1)0^t\,,  &\hbox{ for } a\neq m-1\,,\\[1mm]
 0^s\varphi(p)0^{t+s}\,,           &\hbox{ for } a=m-1\,,
\end{array}
\right.
$$
is, according to (i) of this proposition, also a palindrome, which has a central factor
$(a+1)\varphi(p)0^t(a+1)$ for $a\neq m-1$, and $0\varphi(p)0^t0$ for $a=m-1$.

On the other hand, assume that $b\varphi(p)0^tb\in{\mathcal Pal}(u_\beta)$. If $b\neq 0$,
then using 1. of Remark~\ref{pozn:2.3}, we have
$0^tb\varphi(p)0^tb0^t=\varphi\bigl((b-1)p(b-1)\bigr)0^t\in{\mathcal Pal}(u_\beta)$.
Point (i) implies that $(b-1)p(b-1)\in{\mathcal Pal}(u_\beta)$ and thus $(b-1)p(b-1)$ is
a palindromic extension of $p$. If $b=0$, then Lemma~\ref{l:1} implies that
$10^s\varphi(p)0^t0^s1\in{\mathcal L}(u_\beta)$ and so
$1\varphi\bigl((m-1)p(m-1)0\bigr)\in{\mathcal L}(u_\beta)$, which means that
$(m-1)p(m-1)$ is a palindromic extension of $p$. \pfk

Unlike Arnoux-Rauzy words, in the case of infinite words $u_\beta$ with
$d_\beta(1)=tt\cdots ts$, $t\geq s\geq 2$, it is not difficult to see using
Lemma~\ref{l:1} that there exist palindromes which do not have any palindromic extension.
Such a palindrome is for example the word $0^{t+s-1}$.

\begin{de}
A palindrome $p\in{\mathcal Pal}(u_\beta)$ which has no palindromic extension is called
a {\em maximal} palindrome.
\end{de}

It is obvious that every palindrome is either a central factor of a maximal palindrome, or
is a central factor of palindromes of arbitrary length.

Proposition~\ref{p:1} allows us to define a sequence of maximal palindromes starting from
an initial maximal palindrome. Put
\begin{equation}\label{e:4u}
U^{(1)}\coloneq0^{t+s-1},\qquad U^{(n)}\coloneq\varphi(U^{(n-1)})0^t,\quad\hbox{for }n\geq 2\,.
\end{equation}

Lemma~\ref{l:1} also implies that the palindrome $0^t$ has for $s\geq 2$ two palindromic extensions,
namely $00^t0$ and $10^t1$. Using Proposition~\ref{p:1} we create a sequence of palindromes,
all having two palindromic extensions. Put
\begin{equation}\label{e:4}
V^{(1)}\coloneq0^{t},\qquad V^{(n)}\coloneq\varphi(V^{(n-1)})0^t,\quad\hbox{for }n\geq 2\,.
\end{equation}

\begin{pozn}\label{pozn:citace}
It is necessary to mention that the factors $U^{(n)}$ and $V^{(n)}$ defined above play an
important role in the description of factor complexity of the infinite word $u_\beta$.
Let us cite several results for $u_\beta$ invariant under the substitution~\eqref{e:1}
with $s\geq 2$, taken from~\cite{FrMaPe}, which will be used in the sequel.
\begin{itemize}
\item[{(1)}] Any prefix $w$ of $u_\beta$ is  a left special factor which can be extended
to the left by any letter of the alphabet, i.e.\ $aw\in{\mathcal
L}(u_\beta)$ for all $a\in{\mathcal A}$, or equivalently ${\rm
Lext}(w)=\A$.

\item[{(2)}] Any left special factor $w$ which is not a prefix of $u_\beta$ is a prefix of $U^{(n)}$
for some $n\geq 1$ and such $w$ can be extended to the left by exactly two letters.

\item[{(3)}] The words $U^{(n)}$, $n\geq 1$ are maximal left special factors of $u_\beta$, i.e.\
$U^{(n)}a$ is not a left special factor for any $a\in{\mathcal A}$. The infinite word
$u_\beta$ has no other maximal left special factors.

\item[{(4)}] The word $V^{(n)}$ is the longest common prefix of $u_\beta$ and $U^{(n)}$, moreover,
for every $n\geq 1$ we have
\begin{equation}\label{e:5}
|V^{(n)}|<|U^{(n)}|<|V^{(n+1)}|
\end{equation}

\item[{(5)}] For the first difference of factor complexity we have
$$
\Delta {\mathcal C}(n) = \left\{\begin{array}{cl}
 m &\hbox{ if }\ |V^{(k)}|<n\leq |U^{(k)}|\ \hbox{ for some }k\geq 1\,,\\
 m-1&\hbox{ otherwise}\,.
 \end{array}\right.
$$

\end{itemize}
\end{pozn}

Now we are in position to describe the palindromic extensions in $u_\beta$.
The main result is the following one.

\begin{prop}\label{thm:20}
Let $u_\beta$ be the fixed point of the substitution $\varphi$ given by~\eqref{e:1} with
parameters $t\geq s\geq 2$, and let $p$ be a palindrome in $u_\beta$. Then
\begin{itemize}
\item[(i)] $p$ is a maximal palindrome if and only if $p=U^{(n)}$ for some $n\geq 1$;
\item[(ii)] $p$ has two palindromic extensions in $u_\beta$ if and only if $p=V^{(n)}$
for some $n\geq 1$;
\item[(iii)] $p$ has a unique palindromic extension if and only if $p\neq U^{(n)}$, $p\neq V^{(n)}$
for all $n\geq 1$.
\end{itemize}
\end{prop}

\pf (i) Proposition~\ref{p:1}, point (ii) and the construction of $U^{(n)}$ imply that
 $U^{(n)}$ is a maximal palindrome for every $n$. The proof that no other palindrome $p$ is maximal
 will be done by induction on the length $|p|$ of the palindrome $p$.

 Let $p$ be a maximal palindrome. If $p$ does not contain a non-zero letter, then using
 Lemma~\ref{l:1}, obviously $p=U^{(1)}$. Assume therefore that $p$ contains a non-zero
 letter. Point~1. of Remark~\ref{pozn:2.3} implies that $p=0^t\hat{p}0^t$,
 where $\hat{p}$ is a palindrome. Since $p$ is a maximal palindrome, $\hat{p}$ ends and
 starts in a non-zero letter. Otherwise, $p$ would be extendable to a palindrome, which
 contradicts maximality.
  From 2. of Remark~\ref{pozn:2.3} we obtain that $p=0^t\hat{p}0^t=\varphi(w)0^t$
 for some factor $w$. Proposition~\ref{p:1}, (i), implies that $w$ is a palindrome. Point
 (ii) of the same proposition implies that $w$ has no palindromic extension, i.e.\ $w$ is
 a maximal palindrome, with clearly $|w|<|p|$. The induction hypothesis implies that $w=U^{(n)}$
 for some $n\geq 1$ and $p=\varphi(U^{(n)})0^t=U^{(n+1)}$.

(ii) and (iii) From what we have just proved it follows that every palindrome $p\neq U^{(n)}$,
 $n\geq 1$, has at least one palindromic extension. Since we know that $V^{(n)}$ has
 exactly two palindromic extensions, for proving (ii) and (iii) it remains to show that
 if a palindrome $p$ has more than one extension, then $p=V^{(n)}$, for some $n\geq 1$.

 Assume that $ipi$ and $jpj$ are in
 ${\cal L}(u_\beta)$ for $i,j\in{\cal A}$, $i\neq j$. Obviously, $p$
 is a left special factor of $u_\beta$. We distinguish two cases, according to whether
 $p$ is a prefix of $u_\beta$, or not.

\begin{itemize}
\item
 Let $p$ be a prefix of $u_\beta$. Then there exists a letter $k\in{\cal A}$ such that $pk$
 is a prefix of $u_\beta$ and using (1) of Remark~\ref{pozn:citace}, the word
 $apk\in{\cal L}(u_\beta)$ for every letter $a\in{\cal A}$, in particular $ipk$
 and $jpk$ belong to ${\cal
 L}(u_\beta)$. We have either $k\neq i$, or $k\neq j$; without loss of generality assume
 that $k\neq i$. Since ${\cal L}(u_\beta)$ is closed under reversal, we must have
 $kpi\in{\cal L}(u_\beta)$. Since $ipi$ and $kpi$ are
 in ${\cal L}(u_\beta)$, we obtain that $pi$ is
 also a left special factor of $u_\beta$, and $pi$ is not a prefix of $u_\beta$.
 By (2) of Remark~\ref{pozn:citace}, $p$ is the longest common prefix of $u_\beta$
 and some maximal left special factor $U^{(n)}$, therefore using (4) of Remark~\ref{pozn:citace}
 we have $p=V^{(n)}$.
\item
 If $p$ is a left special factor of $u_\beta$, which is not a prefix of $u_\beta$, then
 by (2) of Remark~\ref{pozn:citace}, $p$ is a prefix of some $U^{(n)}$ and the letters
 $i,j$ are the only possible left extensions of $p$. Since $p\neq U^{(n)}$, there exists
 a unique letter $k$ such that $pk$ is a left special factor of $u_\beta$ and $pk$ is a
 prefix of $U^{(n)}$, i.e.\ the possible left extensions of $pk$ are the letters $i,j$.
 Since by symmetry $kp\in{\cal L}(u_\beta)$, we have $k=i$ or $k=j$, say $k=i$. Since
 $jpk=jpi\in{\cal L}(u_\beta)$, we have also $ipj\in{\cal L}(u_\beta)$. Since by assumption
 $ipi$ and $jpj$ are in ${\mathcal L}(u_\beta)$, both $pi$ and $pj$ are left special factors of $u_\beta$.
 Since $p$ is not a prefix of $u_\beta$, neither $pi$ nor $pj$ are prefixes of $u_\beta$.
 This contradicts the fact that $k$ is a unique letter such that $pk$ is left special.
\end{itemize}

\noindent
 Thus we have shown
 that if a palindrome $p$ has at least two palindromic extensions, then $p=V^{(n)}$.
 \pfk

{}From the above result it follows that if $n\neq |V^{(k)}|$,
$n\neq |U^{(k)}|$ for all $k\geq 1$, then every palindrome of
length $n$ has exactly one palindromic extension, and therefore
${\mathcal P}(n+2)={\mathcal P}(n)$. Inequalities in (4) of
Remark~\ref{pozn:citace} further imply that $|V^{(i)}| \neq |U^{(k)}|$
for all $i,k\geq 1$. Therefore the statement of
Proposition~\ref{thm:20} can be reformulated in the following way:
$$
{\mathcal P}(n+2)-{\mathcal P}(n) = \left\{
\begin{array}{rl}
 1& \hbox{ if }\ n=|V^{(k)}|\,, \\
 -1& \hbox{ if }\ n=|U^{(k)}|\,, \\
 0& \hbox{ otherwise}\,.
\end{array}
 \right.
 $$
Point (5) of Remark~\ref{pozn:citace} can be used for deriving for the second difference
of factor complexity
$$
\Delta^2 {\mathcal C}(n) = \Delta{\mathcal C}(n+1)-\Delta{\mathcal C}(n) = \left\{
\begin{array}{rl}
 1& \hbox{ if }\ n=|V^{(k)}|\,, \\
 -1& \hbox{ if }\ n=|U^{(k)}|\,, \\
 0& \hbox{ otherwise}\,.
\end{array}
\right.
$$
Therefore we have for $s\geq 2$ that ${\mathcal P}(n+2)-{\mathcal
P}(n) = \Delta {\mathcal C}(n+1) - \Delta {\mathcal C}(n)$, for
all $n\in\N$. We thus can derive the following theorem.

\begin{thm}\label{c:sudeliche}
Let $u_\beta$ be the fixed point of the substitution~\eqref{e:1}. Then
$$
{\mathcal P}(n+1)+{\mathcal P}(n) = \Delta{\mathcal C}(n) + 2\,,\qquad\hbox{ for }\
n\in\N\,.
$$
\end{thm}

\pf
 Let the parameter in the substitution~\eqref{e:1} be $s=1$. Then $u_\beta$ is an
 Arnoux-Rauzy word, for which ${\mathcal P}(n+2)-{\mathcal P}(n)=0=
 \Delta{\mathcal C}(n+1)-\Delta{\mathcal C}(n)$. \\
 For $s\geq 2$ we use
 ${\mathcal P}(n+2)-{\mathcal P}(n) = \Delta{\mathcal C}(n+1) - \Delta {\mathcal C}(n)$
 derived above.

 We have
 $$
 \begin{aligned}
{\mathcal P}(n+1)+{\mathcal P}(n) &= {\mathcal P}(0)+{\mathcal P}(1) + \sum_{i=1}^n
\bigl({\mathcal P}(i+1)-{\mathcal P}(i-1)\bigr)  =\\ &=1+m + \sum_{i=1}^n
\bigl(\Delta{\mathcal C}(i)-\Delta{\mathcal C}(i-1)\bigr) = 1+m + \Delta{\mathcal
C}(n)-\Delta{\mathcal C}(0) =\\[3mm]
&= 1+m + \Delta{\mathcal C}(n) - {\mathcal C}(1)+{\mathcal C}(0) = \Delta {\mathcal C}(n)
+ 2\,,
\end{aligned}
 $$
 where we have used ${\mathcal P}(0)={\mathcal C}(0)=1$ and ${\mathcal P}(1)={\mathcal
 C}(1)=m=\#\A$.
 \pfk

\begin{pozn}\label{pozn:omezenost}
According to (5) of Remark~\ref{pozn:citace}, we have $\Delta{\mathcal C}(n)\leq \#\A$.
This implies ${\mathcal P}(n+1)+{\mathcal P}(n)\leq \#\A +2$, and thus the palindromic
complexity is bounded.
\end{pozn}

\section{Centers of palindromes}

We have seen that the set of palindromes of $u_\beta$ is closed under the mapping
$p\mapsto \varphi(p)0^t$. We study the action of this mapping on the centers of the palindromes.
Let us mention that the results of this section are valid for $\beta$ a confluent Pisot number with
$t\geq s\geq1$, i.e.\ also for the Arnoux-Rauzy case.

\begin{de}
Let $p$ be a palindrome of odd length. The {\em center} of $p$ is a letter $a$ such that
$p=wa\widetilde{w}$ for some $w\in{\mathcal A}^*$. The center of a palindrome $p$ of even
length is the empty word.
\end{de}

If palindromes $p_1$, $p_2$ have the same center, then also palindromes $\varphi(p_1)0^t$,
$\varphi(p_2)0^t$ have the same center. This is a consequence of the following lemma.

\begin{lem}\label{l:2}
Let $p,q\in{\mathcal Pal}(u_\beta)$ and let $q$ be a central factor of $p$. Then $\varphi(q)0^t$
is a central factor of $\varphi(p)0^t$.
\end{lem}

Note that the statement is valid also for $q$ being the empty word.

\pf Since $p=wq\widetilde{w}$ for some $w\in{\cal A}^*$, we have $\varphi(p)0^t =
\varphi(w)\varphi(q)\varphi(\widetilde{w})0^t$, which is a palindrome by (i) of
Proposition~\ref{p:1}. It suffices to realize that $0^t$ is a prefix of
$\varphi(\widetilde{w})0^{t}$. Therefore we can write $\varphi(p)0^t =
\varphi(w)\varphi(q)0^t0^{-t}\varphi(\widetilde{w})0^t$. Since
$|\varphi(w)|=|0^{-t}\varphi(\widetilde{w})0^t|$, the word $\varphi(q)0^t$ is a central
factor of $\varphi(p)0^t$.
 \pfk

The following lemma describes the dependence of the center of the palindrome
$\varphi(p)0^t$ on the center of the palindrome $p$. Its proof is a simple application of
properties of the substitution $\varphi$, we will omit it here.

\begin{lem}\label{l:3}
Let $p_1\in{\mathcal Pal}(u_\beta)$ and let $p_2=\varphi(p_1)0^t$.
\begin{itemize}
\item[(i)] If $p_1=w_1a\widetilde{w}_1$, where $a\in{\mathcal A}$, $a\neq m-1$, then
$p_2=w_2(a+1)\widetilde{w}_2$, where $w_2=\varphi(w_1)0^t$.

\item[(ii)] If $p_1=w_1(m-1)\widetilde{w}_1$ and $s+t$ is odd, then
$p_2=w_20\widetilde{w}_2$, where $w_2=\varphi(w_1)0^{\frac{s+t-1}{2}}$.

\item[(iii)] If $p_1=w_1(m-1)\widetilde{w}_1$ and $s+t$ is even, then
$p_2=w_2\widetilde{w}_2$, where $w_2=\varphi(w_1)0^{\frac{s+t}{2}}$.

\item[(iv)] If $p_1=w_1\widetilde{w}_1$ and $t$ is even, then
$p_2=w_2\widetilde{w}_2$, where $w_2=\varphi(w_1)0^{\frac{t}{2}}$.

\item[(v)] If $p_1=w_1\widetilde{w}_1$ and $t$ is odd, then
$p_2=w_20\widetilde{w}_2$, where $w_2=\varphi(w_1)0^{\frac{t-1}{2}}$.
\end{itemize}
\end{lem}

Lemmas~\ref{l:2} and~\ref{l:3} allow us to describe the centers of palindromes $V^{(n)}$
which are in case $s\geq 2$ characterized by having two palindromic extensions.

\begin{prop}\label{p:8}
Let $V^{(n)}$ be palindromes defined by~\eqref{e:4}.
\begin{itemize}
\item[(i)]
If $t$ is even, then for every $n\geq 1$, $V^{(n)}$ has the empty
word $\varepsilon$ for center and $V^{(n)}$ is a central factor of
$V^{(n+1)}$.

\item[(ii)]
If $t$ is odd and $s$ is even, then for every $n\geq 1$, $V^{(n)}$ has
the letter $i\equiv n\!-\!1 \pmod m$ for center, and $V^{(n)}$ is a central factor
of $V^{(n+m)}$.

\item[(iii)]
If $t$ is odd and $s$ is odd, then for every $n\geq 1$, $V^{(n)}$
has the empty word $\varepsilon$ for center if $n\equiv 0
\pmod{(m+1)}$, otherwise it has for center the letter $i\equiv
n\!-\!1\pmod {(m+1)}$. Moreover, $V^{(n)}$ is a central factor of
$V^{(m+n+1)}$.
\end{itemize}
\end{prop}

\pf If $t$ is even, then the empty word $\varepsilon$ is the
center of $V^{(1)}=0^t$. Using Lemma~\ref{l:2} we have that
$\varphi(\varepsilon)0^t=V^{(1)}$ is a central factor of
$\varphi(V^{(1)})0^t=V^{(2)}$. Repeating Lemma~\ref{l:2} we obtain
that $V^{(n)}$ is a central factor of $V^{(n+1)}$. Since
$\varepsilon$ is the center of $V^{(1)}$, it is also the center of
$V^{(n)}$ for all $n\geq 1$.

It $t$ is odd, the palindrome $V^{(1)}$ has center $0$ and using
Lemma~\ref{l:3}, $V^{(2)}$ has center 1, $V^{(3)}$ has center $2$,
\ldots, $V^{(m)}$ has center $m-1$. If moreover $s$ is even, then
$V^{(m+1)}$ has again center $0$. Moreover, from (ii) of
Lemma~\ref{l:3} we see that $0^{s+t}$ is a central factor of
$V^{(m+1)}$, which implies that $V^{(1)}=0^t$ is a central factor
of $V^{(m+1)}$. In case that $s$ is odd, then $V^{(m)}$ having
center $m-1$ implies that $V^{(m+1)}$ has center $\varepsilon$ and
$V^{(m+2)}$ has center $0$. Moreover, using (v) of Lemma~\ref{l:3}
we see that $V^{(1)}=0^t$ is a central factor of $V^{(m+2)}$.
Repeated application of Lemma~\ref{l:2} implies the statement of
the proposition.
 \pfk

As we have said, every palindrome $p$ is either a central factor of a
maximal palindrome $U^{(n)}$, for some $n\geq 1$, or $p$ is a
central factor of palindromes with increasing length. An example
of such a palindrome is $V^{(n)}$, for $n\geq 1$, which is according to Proposition~\ref{p:8}
central factor of palindromes of arbitrary length.
According to the notation introduced by Cassaigne in~\cite{cassaigne}
for left and right special factors extendable to arbitrary length special factors, we
introduce the notion of infinite palindromic branch. We will study infinite palindromic branches
in the next section.

\section{Infinite palindromic branches}

\begin{de}
Let $v=\cdots v_3v_2v_1$ be a left infinite word in the alphabet ${\mathcal A}$.
Denote by $\widetilde{v}$ the right infinite word $\widetilde{v}=v_1v_2v_3\cdots$.
\begin{itemize}
\item Let $a\in{\mathcal A}$.
If for every index $n\geq 1$, the word $p=v_nv_{n-1}\cdots v_1av_1v_2\cdots v_n\in{\mathcal
Pal}(u_\beta)$, then the two-sided infinite word $va\widetilde{v}$ is called an
infinite palindromic branch of $u_\beta$ with center $a$, and the palindrome $p$ is
called a central factor of the infinite palindromic branch $va\widetilde{v}$.

\item If for every index $n\geq 1$, the word $p=v_nv_{n-1}\cdots v_1v_1v_2\cdots v_n\in{\mathcal Pal}(u_\beta)$,
then the two-sided infinite word $v\widetilde{v}$ is called an
infinite palindromic branch of $u_\beta$ with center
$\varepsilon$, and the palindrome $p$ is called a central factor
of the infinite palindromic branch $v\widetilde{v}$.
\end{itemize}
\end{de}

Since for Arnoux-Rauzy words every palindrome has exactly one
palindromic extension, we obtain for every letter $a\in{\mathcal
A}$ exactly one infinite palindromic branch with center $a$; there
is also one infinite palindromic branch with center $\varepsilon$.

Obviously, every infinite word with bounded palindromic complexity ${\mathcal P}(n)$ has only
a finite number of infinite palindromic branches. This is therefore valid also for $u_\beta$.

\begin{prop}\label{p:9}
The infinite word $u_\beta$ invariant under the
substitution~\eqref{e:1} has for each center $c\in{\mathcal
A}\cup\{\varepsilon\}$ at most one infinite palindromic branch
with center $c$.
\end{prop}

\pf Lemma~\ref{l:3} allows us to create from one infinite
palindromic branch another infinite palindromic branch. For
example, if $va\widetilde{v}$ is an infinite palindromic branch
with center $a\neq m-1$, then using (i) of Lemma~\ref{l:3}, the
two-sided word $\varphi(v)0^t(a+1)0^t\widetilde{\varphi(v)}$ is an
infinite palindromic branch with center $(a+1)$. Similarly for the
center $m-1$ or $\varepsilon$. Obviously, this procedure creates
from distinct palindromic branches with the same center $c\in{\cal
A}\cup\{\varepsilon\}$ again distinct palindromic branches, for
which the length of the maximal common central factor is longer
than the length of the maximal common central factor of the
original infinite palindromic branches. This would imply that
$u_\beta$ has infinitely many infinite palindromic branches, which
is in contradiction with the boundedness of the palindromic
complexity of $u_\beta$, see~Remark~\ref{pozn:omezenost}.
 \pfk

\begin{pozn}\label{pozn:infpalbr}
Examples of infinite palindromic branches can be easily obtained from
Proposition~\ref{p:8} as a two-sided limit of palindromes $V^{(k_n)}$ for a suitably
chosen subsequence $(k_n)_{n\in\N}$ and $n$ going to infinity, namely
\begin{itemize}
\item
If $t$ is even, then the two-sided limit of palindromes $V^{(n)}$
is an infinite palindromic branch with center $\varepsilon$.

\item
If $t$ is odd and $s$ even, then the two-sided limit of palindromes $V^{(k+mn)}$ for
$k=1,2,\ldots,m$ is an infinite palindromic branch with center $k-1$.

\item
If $t$ is odd and $s$ odd, then the two-sided limit of palindromes
$V^{(k+(m+1)n)}$ for $k=1,2,\ldots,m$ is an infinite palindromic
branch with center $k-1$, and for $k=m+1$ it is an infinite
palindromic branch with center $\varepsilon$.
\end{itemize}
\end{pozn}

\begin{coro}\label{c}\
\begin{itemize}
\item[(i)] If $s$ is odd, then $u_\beta$ has exactly one infinite palindromic branch with center $c$
for every $c\in{\mathcal A}\cup\{\varepsilon\}$.
\item[(ii)] If $s$ is even and $t$ is odd, then $u_\beta$ has exactly one infinite palindromic branch
with center $c$ for every $c\in{\mathcal A}$, and $u_\beta$ has no
infinite palindromic branch with center $\varepsilon$.
\item[(iii)] If $s$ is even and $t$ is even, then $u_\beta$ has exactly one infinite palindromic branch
with center $\varepsilon$, and $u_\beta$ has no infinite
palindromic branch with center $a\in{\mathcal A}$.
\end{itemize}
\end{coro}

\pf
 According to Proposition~\ref{p:9}, $u_\beta$ may have at most one infinite palindromic
 branch for each center $c\in{\mathcal A}\cup\{\varepsilon\}$. Therefore it suffices to show
 existence/non-existence of such a palindromic branch. We distinguish four cases:

--- Let $s$ be odd and $t$ odd. Then an infinite palindromic branch with center $c$
 exists for every $c\in{\mathcal A}\cup\{\varepsilon\}$, by Remark~\ref{pozn:infpalbr}.

--- Let $s$ be odd and $t$ even. The existence of an infinite
palindromic branch with center $\varepsilon$ is ensured again by
Remark~\ref{pozn:infpalbr}. For determining the infinite
palindromic branches with other centers, we define a sequence of
words
$$
W^{(1)} = 0\,,\qquad\qquad W^{(n+1)} = \varphi(W^{(n)})0^t\,,\quad n\in\N, n\geq 1\,.
$$
Since $s+t$ is odd, using (i) and (ii) of Lemma~\ref{l:3}, we know that $W^{(n)}$ is a
palindrome with  center $i\equiv n\!-\!1 \pmod m$. In particular, we
have that $0=W^{(1)}$
is a central factor of $W^{(m+1)}$. Using Lemma~\ref{l:2}, also $W^{(n)}$ is a central
factor of $W^{(m+n)}$ for all $n\geq 1$. Therefore we can construct
 the two-sided limit of
palindromes $W^{(k+mn)}$ for $n$ going to infinity, to obtain an infinite palindromic
branch with center $k-1$ for all $k=1,2,\ldots,m$.

--- Let $s$ be even and $t$ be odd. Then an infinite palindromic branch with center $c$
 exists for every $c\in{\mathcal A}$, by Remark~\ref{pozn:infpalbr}.
 A palindromic branch with center $\varepsilon$ does not exist, since using
 Lemma~\ref{l:1} two non-zero letters in the
 word $u_\beta$ are separated by a block of $0$'s of odd length, which implies that
 palindromes of even length must be shorter than $t+s$.

--- Let $s$ and $t$ be even. The existence of an infinite
palindromic branch with center $\varepsilon$ is ensured again by
Remark~\ref{pozn:infpalbr}. Infinite palindromic
 branches with other centers do not exist. The reason is that in this case the maximal
 palindrome $U^{(1)}=0^{t+s-1}$ has center
 $0$ and using Lemma~\ref{l:3} the palindromes $U^{(2)}$, $U^{(3)}$, \ldots, $U^{(m)}$
 have centers $1,2,\ldots,m-1$, respectively. For all $n>m$ the center of $U^{(n)}$ is
 the empty word $\varepsilon$. If there existed an infinite palindromic branch
 $va\widetilde{v}$, then the maximal common central factor $p$ of $va\widetilde{v}$ and $U^{(a+1)}$
 would be a palindrome with center $a$ and with two palindromic extensions. Using
 Proposition~\ref{thm:20}, $p=V^{(k)}$ for some $k$. Proposition~\ref{p:8} however implies
 that for $t$ even the center of $V^{(k)}$ is the empty word $\varepsilon$, which is a
 contradiction.
\pfk

\begin{pozn}\label{pozn:23}
The proof of the previous corollary implies:

\begin{itemize}
\item[(i)]
In case $t$ odd, $s$ even, $u_\beta$ has only finitely many palindromes of even length,
all of them being central factors of $U^{(1)}=0^{t+s-1}$.

\item[(ii)]
In case $t$ and $s$ are
even, $u_\beta$ has only finitely many palindromes of odd length
and all of them are central factors of one of the palindromes $U^{(1)}$, $U^{(2)}$,
\ldots, $U^{(m)}$, with center $0,1,\ldots,m-1$, respectively.
\end{itemize}
\end{pozn}

\section{Palindromic complexity of $u_\beta$}

The aim of this section is to give explicit values of the palindromic complexity of
$u_\beta$. We shall derive them from Theorem~\ref{c:sudeliche}, which expresses
${\mathcal P}(n)+ {\mathcal P}(n+1) $ using the first difference of factor complexity;
and from (5) of Remark~\ref{pozn:citace}, which recalls the results about ${\mathcal
C}(n)$ of~\cite{FrMaPe}.

\begin{thm}\label{thm:10}
Let $u_\beta$ be the fixed point of the substitution~\eqref{e:1}, with parameters $t\geq s\geq2$.
\begin{itemize}
\item[(i)] Let $s$ be odd and let $t$ be even. Then
\begin{eqnarray*}
{\mathcal P}(2n+1)&=&m\\[1mm]
{\mathcal P}(2n)&=&\left\{
\begin{array}{cl}
2,&\hbox{ if } \ |V^{(k)}|<2n\leq |U^{(k)}|\, \hbox{ for some } k,\\[2mm]
1,&\hbox{ otherwise}.
\end{array}\right.
\end{eqnarray*}

\item[(ii)] Let $s$ and $t$ be odd. Then
\begin{eqnarray*}
{\mathcal P}(2n+1)&=&\left\{
\begin{array}{cl}
m+1,&\hbox{ if } \ |V^{(k)}|<2n+1\leq |U^{(k)}|\, \hbox{ for some $k$ } \\
& \hbox{ with } k\not\equiv 0 \pmod{(m+1)},\\[2mm]
m,&\hbox{ otherwise}.
\end{array}\right.\\[3mm]
{\mathcal P}(2n)&=&\left\{
\begin{array}{cl}
2,&\hbox{ if } \ |V^{(k)}|<2n\leq |U^{(k)}|\, \hbox{ for some $k$ }\\
& \hbox{ with } k\equiv 0 \pmod{(m+1)},\\[2mm]
1,&\hbox{ otherwise}.
\end{array}\right.
\end{eqnarray*}

\item[(iii)] Let $s$ be even and $t$ be odd. Then
\begin{eqnarray*}
{\mathcal P}(2n+1)&=&\left\{
\begin{array}{cl}
m+2,&\hbox{ if }\ |V^{(k)}|<2n+1\leq |U^{(k)}|\, \hbox{ for some } k\geq 2\,,\\[2mm]
m,&\hbox{ if } \ 2n+1\leq |V^{(1)}|\,,\\[2mm]
m+1,&\hbox{ otherwise}.
\end{array}\right.\\[3mm]
{\mathcal P}(2n)&=&\left\{
\begin{array}{cl}
1,&\hbox{ if }\ 2n\leq|U^{(1)}|\,,\\[2mm]
0,&\hbox{ otherwise}.\\

\end{array}\right.
\end{eqnarray*}

\item[(iv)] Let $s$ and $t$ be even. Then
\begin{eqnarray*}
{\mathcal P}(2n+1)&=&\left\{
\begin{array}{cl}
\#\bigl\{k\leq m\,\bigm|\, 2n+1\leq |U^{(k)}|\bigr\},&\hbox{ if } \ 2n+1\leq
|U^{(m)}|\,,\\[2mm]
0,&\hbox{ otherwise}.
\end{array}\right.\\[3mm]
{\mathcal P}(2n)&=&\left\{
\begin{array}{cl}
m+2,&\hbox{ if } \ |V^{(k)}|<2n\leq |U^{(k)}|\, \\
&\hbox{ for some } k\geq m+1\,,\\[2mm]
\#\bigl\{k\leq m\,\bigm|\, 2n > |V^{(k)}|\bigr\}+1,&\hbox{ if } \ 2n\leq |V^{(m+1)}|\,,\\[2mm]
m+1,&\hbox{ otherwise}.
\end{array}\right.
\end{eqnarray*}

\end{itemize}
\end{thm}

\pf We prove the statement by cases:

\begin{itemize}
\item[(i)] Let $s$ be odd and $t$ be even. It is enough to show that ${\mathcal P}(2n+1)=m$ for
all $n\in\N$. The value of ${\mathcal P}(2n)$ can then be easily calculated from
Theorem~\ref{c:sudeliche} and (5) of Remark~\ref{pozn:citace}.

{}From (i) of Corollary~\ref{c} we know that there exists an
infinite palindromic branch with center $c$ for all $c\in\A$. This
implies that ${\mathcal P}(2n+1)\geq m$. In order to show the
equality, it suffices to show that all maximal palindromes
$U^{(k)}$ are of even length, or equivalently, have $\varepsilon$
for center. Since both $t$ and $t+s-1$ are even, $0^t=V^{(1)}$ is
a central factor of $0^{t+s-1}=U^{(1)}$. Using Lemma~\ref{l:2},
$V^{(k)}$ is a central factor of $U^{(k)}$ for all $k\geq 1$.
According to (i) of Proposition~\ref{p:8}, $V^{(k)}$ are
palindromes of even length, and thus also the maximal palindromes
$U^{(k)}$ are of even length. Therefore they do not contribute to
${\mathcal P}(2n+1)$.

\item[(ii)]  Let $s$ and $t$ be odd. We shall determine ${\mathcal P}(2n)$ and the values of
${\mathcal P}(2n+1)$ can be deduced from Theorem~\ref{c:sudeliche} and (5) of
Remark~\ref{pozn:citace}.

{}From (i) of Corollary~\ref{c} we know that there exists an
infinite palindromic branch with center $\varepsilon$. Thus
${\mathcal P}(2n)\geq 1$ for all $n\in\N$. Again, $V^{(1)}=0^t$ is
a central factor of $U^{(1)}=0^{t+s-1}$, and thus $V^{(k)}$ is a
central factor of $U^{(k)}$ for all $k\geq 1$. A palindrome of even
length, which is not a central factor of an infinite palindromic
branch must be a central factor of $U^{(k)}$ for some $k$, and
longer than $|V^{(k)}|$. Since $|U^{(k)}|<|V^{(k+1)}|<|U^{(k+1)}|$
(cf. (5) of Remark~\ref{pozn:citace}), at most one such palindrome
exists for each length. We have ${\mathcal P}(2n)\leq 2$. It
suffices to determine for which $k$, the maximal palindrome
$U^{(k)}$ is of even length, which happens exactly when its
central factor $V^{(k)}$ is of even length and that is, using
(iii) of Proposition~\ref{p:8}, for $k\equiv 0\pmod{(m+1)}$.

\item[(iii)] Let $s$ be even and $t$ be odd. According to (i) of Remark~\ref{pozn:23}, all
palindromes of even length are central factors of $U^{(1)}=0^{t+s-1}$. Therefore
${\mathcal P}(2n)=1$ if $2n\leq |U^{(1)}|$ and 0 otherwise. The value of ${\mathcal
P}(2n+1)$ can be calculated from Theorem~\ref{c:sudeliche} and (5) of
Remark~\ref{pozn:citace}.

\item[(iv)] Let $s$ and $t$ be even. Using (ii) of Remark~\ref{pozn:23}, the only
palindromes of odd length are central factors of $U^{(k)}$ for $k=1,2,\ldots,m$. Therefore
${\mathcal P}(2n+1)=0$ for $2n+1>|U^{(m)}|$. If $2n+1\leq|U^{(m)}|$, the number of
palindromes of odd length is equal to the number of maximal palindromes longer than
$2n+1$. The value of ${\mathcal P}(2n)$ can be calculated from
Theorem~\ref{c:sudeliche} and (5) of Remark~\ref{pozn:citace}.
\end{itemize}
\pfk

For the determination of the value ${\mathcal P}(n)$ for a given $n$, we have to know
$|V^{(k)}|$, $|U^{(k)}|$. In~\cite{FrMaPe} it is shown that
$$
|V^{(k)}|=t\sum_{i=0}^{k-1}G_i\,, \qquad\hbox{ and }\qquad
|U^{(k)}|=|V^{(k)}|+(s-1)G_{k-1}\,,
$$
where $G_n$ is a sequence of integers defined by the recurrence
$$
\begin{array}{ll}
G_0=1\,,\quad G_n=t(G_{n-1}+\cdots+G_0)+1\,, &\hbox{ for } 1\leq n\leq m-1\,,\\[2mm]
G_n=t(G_{n-1}+\cdots+G_{n-m+1}) + sG_{n-m}\,, &\hbox{ for } n\geq m\,.
\end{array}
$$
The sequence $(G_n)_{n\in\N}$ defines the canonical linear numeration system
associated with the number $\beta$, see~\cite{bertrand2} for general results
on these numeration systems. In this particular case, $(G_n)_{n\in\N}$ defines
a confluent linear numeration system, see~\cite{frougny} for its properties.

\section{Substitution invariance of palindromic branches}

Infinite words $u_\beta$ are invariant under the
substitution~\eqref{e:1}. One can ask whether also their infinite
palindromic branches are invariant under a substitution. In case
that an infinite palindromic branch has as its center the empty
word $\varepsilon$, we can use the notion of invariance under
substitution as defined for pointed two-sided infinite words. We
restrict our attention to infinite palindromic branches of such
type.

Recall that an infinite palindromic branch of $u_{\beta}$ with center $\varepsilon$ exists,
(according to Corollary~\ref{c}), only if in the R\'enyi
expansion $d_\beta(1)=tt\cdots ts$, $t$ is even, or both $t$ and $s$ are
odd. Therefore we shall study only such parameters.

Let us first study the most simple case, $d_\beta(1)=t1$ for
$t\geq 1$. Here $\beta$ is a quadratic unit, and the infinite word
$u_\beta$ is a Sturmian word, expressible in the form of the
mechanical word $\mu_{\alpha,\varrho}$,
$$
\mu_{\alpha,\varrho}(n)=\bigl\lfloor(n+1)\alpha+\varrho\bigr\rfloor
- \bigl\lfloor n\alpha+\varrho\bigr\rfloor\,,\quad n\in\N\,,
$$
where the irrational slope $\alpha$ and the intercept $\rho$
satisfy $\alpha=\rho=\frac{\beta}{\beta+1}$. The infinite
palindromic branch with center $\varepsilon$ of the above word
$u_\beta=\mu_{\alpha,\varrho}$ is a two-sided Sturmian word with
the same slope $\alpha=\frac{\beta}{\beta+1}$, but
intercept~$\frac12$. Indeed, two mechanical words with the same
slope  have the same set of factors independently on their
intercepts, and moreover the sturmian word $\mu_{\alpha,\frac12}$
is an infinite palindromic branch of itself, since
$$
\mu_{\alpha,\frac12}(n)=\mu_{\alpha,\frac12}(-n-1)\,,\quad\hbox{
for all } n\in\Z\,.
$$
Therefore if
$v=\mu_{\alpha,\frac12}(0)\mu_{\alpha,\frac12}(1)\mu_{\alpha,\frac12}(2)\cdots$,
then  $\widetilde{v}v$  is the  infinite palindromic branch of
$u_\beta$ with the center $\varepsilon$.

Since the Sturmian word $\mu_{\alpha,\varrho}$ coincides with
$u_\beta$, it is invariant under the substi\-tution $\varphi$. As
a consequence of~\cite{parvaix}, the slope $\alpha$ is a Sturm
number, i.e.\ a quadratic number in $(0,1)$ such that its
conjugate $\alpha'$ satisfies $\alpha'\notin(0,1)$, (using the
equivalent definition of Sturm numbers given in~\cite{allauzen}).

The question about the substitution invariance of the infinite
palindromic branch $\widetilde{v}v$ is answered using the result
of~\cite{BaMaPe} (or also~\cite{yasutomi,berthe}). It says that a
Sturmian word whose slope is a Sturm number, and whose intercept
is equal to $\frac12$, is substitution invariant as a two-sided
pointed word, i.e.\ there exists a substitution $\psi$ such that
$\widetilde{v}|v=\psi(\widetilde{v})|\psi(v)$.

\begin{ex}
The Fibonacci word $u_\beta$ for $d_\beta(1)=11$ is a fixed point of the substitution
$$
\varphi(0)=01,\quad \varphi(1)=0\,.
$$
Its infinite palindromic branch with center $\varepsilon$ is
$$
\widetilde{v}v\qquad\hbox{ for }\qquad v=010100100101001001010\cdots
$$
which is the fixed point $\lim_{n\to\infty}\psi^n(0)|\psi^n(0)$ of the substitution
$$
\psi(0)=01010,\quad \psi(1)=010\,.
$$
\end{ex}

\medskip

Let us now study the question whether infinite palindromic
branches in $u_\beta$ for general $d_\beta(1)=tt\cdots ts$ with
$t$ even, or $t$ and $s$ odd, are also substitution invariant. It
turns out that the answer is positive. For construction of a
substitution $\psi$ under which a given palindromic branch is
invariant, we need the following lemma.

\begin{lem}\label{l:24}
Let $v\widetilde{v}$ be an infinite palindromic branch with center
$\varepsilon$. Then the left infinite word $v=\cdots v_3v_2v_1$
satisfies
$$
\begin{array}{ll}
v=\varphi(v)0^{\frac{t}2} & \hbox{ for $t$ even,}\\[1mm]
v=\varphi^{m+1}(v)\varphi^{m}(0^{\frac{t+1}{2}})0^{\frac{t-s}{2}} & \hbox{ for $t$ and
$s$ odd.}
\end{array}
$$
\end{lem}

\pf Let $t$ be even and let $v\widetilde{v}$ be the unique
infinite palindromic branch with center $\varepsilon$. Recall that
$v\widetilde{v}=\lim_{n\rightarrow\infty}V^{(n)}$.
Consider arbitrary suffix $v_{\rm suf}$ of $v$, i.e.\ $v_{\rm
suf}\widetilde{v}_{\rm suf}$ is a palindrome of $u_\beta$ with
center $\varepsilon$. Denote $w\coloneq\varphi(v_{\rm
suf})0^{\frac{t}{2}}$. Using (iv) of Lemma~\ref{l:3} the word
$p=w\widetilde{w}$ is a palindrome of $u_\beta$ with center
$\varepsilon$. We show by contradiction that $w$ is a suffix of $v$.

Suppose that $p = w\widetilde{w}$ is not a central factor of
$\lim_{n\rightarrow\infty}V^{(n)}$, then there exists a unique $n$
such that $p$ is a central factor of $U^{(n)}$. Then according to
Proposition \ref{thm:20}, $p$ is uniquely extendable into a
maximal palindrome. In that case we take a longer suffix $v'_{\rm
suf}$ of $v$, so that the length of the palindrome
$p'=w'\widetilde{w}'$, $w'\coloneq\varphi(v'_{\rm
suf})0^{\frac{t}{2}}$ satisfies $|p'|>|U^{(n)}|$. However, $p'$
(since it contains $p$ as its central factor) is a palindromic
extension of $p$, and therefore $p'$ is a central factor of
$U^{(n)}$, which is a contradiction. Thus $\varphi(v_{\rm
suf})0^{\frac{t}{2}}$ is a suffix of $v$ for all suffixes $v_{\rm
suf}$ of $v$, therefore $v=\varphi(v)0^{\frac{t}{2}}$.

Let now $s$ and $t$ be odd. If $v_{\rm suf}$ is a suffix of the
word $v$, then $v_{\rm suf}\widetilde{v}_{\rm suf}$ is a
palindrome of $u_\beta$ with center $\varepsilon$. Using
Lemma~\ref{l:3}, the following holds true.
$$
\begin{array}{rclccl}
w_0&=&\varphi(v_{\rm suf})0^{\frac{t-1}{2}} &\qquad\implies\qquad&
w_00\widetilde{w}_0 &\in {\mathcal Pal}(u_\beta)\\[1mm]
w_1&=&\varphi(w_0)0^{t} &\qquad\implies\qquad&
w_11\widetilde{w}_1 &\in {\mathcal Pal}(u_\beta)\\[1mm]
w_2&=&\varphi(w_1)0^{t} &\qquad\implies\qquad&
w_22\widetilde{w}_2 &\in {\mathcal Pal}(u_\beta)\\[1mm]
&&&\qquad\vdots\qquad&\\[1mm]
w_{m-1}&=&\varphi(w_{m-2})0^{t} &\qquad\implies\qquad&
w_{m-1}(m-1)\widetilde{w}_{m-1} &\in {\mathcal Pal}(u_\beta)\\[1mm]
w_{\varepsilon}&=&\varphi(w_{m-1})0^{\frac{s+t}{2}}
&\qquad\implies\qquad& w_{\varepsilon}\widetilde{w}_{\varepsilon}
&\in {\mathcal Pal}(u_\beta)
\end{array}
$$
Together we obtain
$$
w_\varepsilon = \varphi^{m+1}(v_{\rm
suf})\varphi^{m}(0^{\frac{t-1}{2}})\varphi^{m-1}(0^t)\cdots
\varphi^2(0^t)\varphi(0^t)0^{\frac{s+t}{2}}\,.
$$
Since $\varphi^m(0)= \varphi^{m-1}(0^{t})\varphi^{m-2}(0^t)\cdots
\varphi(0^t)0^{s}$, the word $w_\varepsilon$ can be rewritten in a
simpler form
$$
w_\varepsilon= \varphi^{m+1}(v_{\rm suf})
\varphi^{m}(0^{\frac{t-1}{2}})\varphi^{m}(0)0^{\frac{t-s}{2}} =
\varphi^{m+1}(v_{\rm suf})
\varphi^{m}(0^{\frac{t+1}{2}})0^{\frac{t-s}{2}}
$$
Since $w_\varepsilon$ is again a suffix of $v$, the statement of
the lemma for $s$ and $t$ odd holds true. \pfk

\begin{thm}\label{t:subst}
Let $u_\beta$ be the fixed point of the substitution $\varphi$
given by~\eqref{e:1}, and let $v\widetilde{v}$ be the infinite
palindromic branch of $u_\beta$ with center $\varepsilon$. Then
the left-sided infinite word $v$ is invariant under the
substitution $\psi$ defined for all letters
$a\in\{0,1,\ldots,m-1\}$ by
$$
\psi(a) =\left\{\begin{array}{lll}
w^{-1}\varphi(a)w\,,& \hbox{ where }\ w=0^{\frac{t}{2}}\,,&\hbox{ for $t$ even,}\\[2mm]
w^{-1}\varphi^{m+1}(a)w\,,& \hbox{ where }\
w=\varphi^m(0^{\frac{t+1}2})0^{\frac{t-s}2}\,, &\hbox{ for $t$ and
$s$ odd.}
\end{array}\right.
$$
Moreover, if $t$ is even, then $\psi(a)$ is a palindrome for all
$a\in \A$ and $v\widetilde{v}$
 as a pointed sequence is invariant under the same substitution
 $\psi$.
\end{thm}

\pf

First let us show that the substitution $\psi$ is well defined.

\begin{itemize}
\item Let $t$ be even. Since $0^{\frac{t}2}$ is a prefix of
$\varphi(a)$ for all $a\in\{0,1,\ldots,m-2\}$ and $\varphi(m-1)=0^s$, therefore $0^{\frac{t}2}$
is a prefix of $\varphi(m-1)0^{\frac{t}2}=0^{s+\frac{t}2}$.

\item Let $t$ and $s$ be odd. Let us verify that $w$ is a prefix of $\varphi^{m+1}(a)w$.

-- If $a\neq m-1$, we show that $w=\varphi^m(0^{\frac{t+1}2})0^{\frac{t-s}2}$ is a prefix of
$$
\varphi^{m+1}(a) = \varphi^m\bigl(0^t(a+1)\bigr) =
\varphi^m(0^{\frac{t+1}2}) \varphi^m(0^{\frac{t-1}2})
\varphi^m(a+1)\,.
$$
It suffices to show that $0^{\frac{t-s}2}$ is a prefix of
$\varphi^m(0^{\frac{t-1}2})$. For $t=s$ it is obvious. For
$t>s\geq 1$ we obtain $t\geq 3$ and so
$\varphi^m(0^{\frac{t-1}2})=\varphi^m(0)\varphi^m(0^{\frac{t-3}2})$
and clearly $0^{\frac{t-s}2}$ is a prefix of $\varphi^m(0)$.

-- If $a=m-1$, then
$$
\varphi^{m+1}(m-1)w = \varphi^m(0^s)\varphi^m(0^{\frac{t+1}2})0^{\frac{t-s}{2}} =
\varphi^m(0^{\frac{t+1}2})\varphi^m(0)\varphi^m(0^{s-1})0^{\frac{t-s}{2}}\,.
$$
Since $0^{\frac{t-s}{2}}$ is a prefix of $\varphi^m(0)$, the correctness of the definition
of the substitution $\psi$ is proven.
\end{itemize}

\noindent
Now it is enough to prove that $\psi(v)=v$.
Lemma~\ref{l:24} says that in the case that $t$ is even the left infinite word
$v=\cdots v_{3}v_{2}v_{1}$ satisfies $v=\varphi(v)w$. Thus we have
$$
\begin{aligned}
\psi(v)=\cdots \psi(v_3)\psi(v_2)\psi(v_1) &=
\cdots w^{-1} \varphi(v_3)ww^{-1}\varphi(v_2) ww^{-1}\varphi(v_2) w =\\
&=\cdots \varphi(v_3)\varphi(v_2)\varphi(v_1)w
\ = \ \varphi(v) w \ = \ v\,.
\end{aligned}
$$
In case that $t$ and $s$ are odd, the proof is the same, using $\varphi^{m+1}$ instead of $\varphi$.

If $t$ is even, it is clear from the prescription for $\psi$, that
$\psi(a)$ is a palindrome for any letter $a$, which implies the
invariance of the word $v\widetilde{v}$ under  $\psi$. \pfk

Let us mention that  for $t,s$ odd the words $\psi(a)$, $a\in \A$,
may not be palindromes. In that case the right-sided word
$\widetilde{v}$ is invariant under another substitution, namely
$a\mapsto \widetilde{\psi(a)}$. Nevertheless even for $t,s$ odd it
may  happen that $\psi(a)$ is a palindrome for all letters. Then
the two-sided word $v\widetilde{v}$ is invariant under $\psi$.
This situation is illustrated  on the following example.

\begin{ex}\label{ex:?}
Consider the Tribonacci word, i.e.\ the word $u_\beta$ for $d_\beta(1)=111$. It is the fixed point
of the substitution
$$
\varphi(0)=01,\quad \varphi(1)=02,\quad \varphi(2)=0\,,
$$
which is in the form~\eqref{e:1} for $t=s=1$ and $m=3$. Therefore
$w=\varphi^3(0)=0102010$.
The substitution $\psi$, under which the infinite palindromic branch
$v\widetilde{v}$ of the Tribonacci word is invariant, is therefore given as
$$
\begin{array}{rclcl}
\psi(0)&\coloneq& w^{-1}\varphi^4(0)w &=&0102010102010,\\
\psi(1)&\coloneq& w^{-1}\varphi^4(1)w &=&01020102010,\\
\psi(2)&\coloneq& w^{-1}\varphi^4(2)w &=&0102010.
\end{array}
$$
Note that the substitution $\psi$ has
the following property: the word $\psi(a)$ is a palindrome for every $a\in\A$.
\end{ex}

\section{Number of palindromes in the prefixes of $u_\beta$}

In~\cite{DrJuPi} the authors obtain an interesting result which
says that every finite word $w$
contains at most $|w|+1$ different palindromes. (The empty word is considered as a
palindrome contained in every word.) Denote by $P(w)$ the number of palindromes contained
in the finite word $w$. Formally, we have
$$
P(w) \leq |w|+1\qquad \hbox{ for every finite word } w.
$$
The finite words $w$ for which the equality is reached are called {\em full} (as suggested
in~\cite{brlek}). An infinite word $u$ is called full, if all its prefixes are full.
In~\cite{DrJuPi} the authors have shown that every Sturmian word is full. They have shown the same
property for episturmian words.

The infinite word $u_\beta$ can be full only if its language is closed under reversal,
i.e.\ in the simple Parry case for $d_\beta(1)=tt \cdots ts$, $t\geq s\geq 1$. For $s\geq
2$ such words are not episturmian, nevertheless, we shall show that they are full.

We shall use the notions and results introduced in~\cite{DrJuPi}.

\begin{de}
A finite word $w$ satisfies property ${\it Ju}$, if there exists a palindromic suffix of $w$ which is
unioccurrent in $w$.
\end{de}

Clearly, if $w$ satisfies ${\it Ju}$, then it has exactly one palindromic suffix which is
unioccurrent, namely the longest palindromic suffix of $w$.

\begin{prop}[\cite{DrJuPi}]\label{p:drjupi}
Let $w$ be a finite word. Then $P(w)=|w|+1$ if and only if all the prefixes $\hat{w}$
of $w$
satisfy ${\it Ju}$, i.e.\ have a palindrome suffix which is unioccurrent in $\hat{w}$.
\end{prop}

\begin{thm}\label{thm:25}
The infinite word $u_\beta$ invariant under the substitution~\eqref{e:1} is full.
\end{thm}

\pf We show the statement using Proposition~\ref{p:drjupi} by
contradiction. Let $w$ be a prefix of $u_\beta$ of minimal length
which does not satisfy ${\it Ju}$, and let $X0^k$ is a suffix of
$w$ with $X\neq0$.

First we show that $k\in\{0,t+1\}$. For, if $1\leq k\leq t$ or
$t+2\leq k$, then $q$ is the maximal palindromic suffix of
$w0^{-1}$ if and only if $0q0$ is the maximal palindromic suffix
of $w$. Since $0q0$ occurs at least twice in $w$, then also $q$
occurs at least twice in $w0^{-1}$, which is a contradiction with
the minimality of $w$.

Define
$$
w_1=\left\{
\begin{array}{ll}
w0^{t} &\hbox{ if $w$ has suffix } X\neq0\,, \\[2mm]
w0^{s-1} &\hbox{ if $w$ has suffix } X0^{t+1},\ X\neq0\,.
\end{array}\right.
$$
For the maximal palindromic suffix $p$ of $w$ denote
$$
p_1=\left\{
\begin{array}{ll}
0^{t}p0^t &\hbox{ if $w$ has suffix } X\neq0\,, \\[2mm]
0^{s-1}p0^{s-1} &\hbox{ if $w$ has suffix } X0^{t+1},\ X\neq0\,.
\end{array}\right.
$$

Since in $u_\beta$ every two non-zero letters are separated by the
word $0^t$ or $0^{t+s}$, we obtain that
\begin{itemize}
\item[(i)]
 $p_1$ is the maximal palindromic suffix  of $w_1$.
\item[(ii)]
 the position of centers of palindromes $p$ and $p_1$ coincide
 in all occurrences in $u_\beta$.
\end{itemize}
Since $p$ occurs in $w$ at least twice, also the palindromic
suffix $p_1$ occurs at least twice in $w_1$, i.e.\  the word $w_1$
is a prefix of $u_\beta$ which does not satisfy ${\it Ju}$.

{}From the definition of $w_1$ it follows that
$$
w_1=\varphi(\hat{w})0^t
$$
for some prefix $\hat{w}$ of $u_\beta$. Thus the maximal
palindromic suffix $p_1$ of $w_1$ is of the form
$p_1=\varphi(\hat{p})0^t$, where $\hat{p}$ is a factor of
$\hat{w}$. According to (i) of Proposition~\ref{p:1}, $\hat{p}$ is
a palindrome, and the same proposition implies that $\hat{p}$ is
the maximal palindromic suffix of $\hat{w}$. Since $p_1$ occurs at
least twice in $w_1$, also $\hat{p}$ occurs at least twice in
$\hat{w}$. Therefore $\hat{w}$ does not satisfy the property ${\it
Ju}$. As
$$
|\hat{w}|<|\varphi(\hat{w})|<|w|,
$$
we have a contradiction with the minimality of $w$.
\pfk

\section{Conclusions}

The study of palindromic complexity of an uniformly recurrent
infinite word is interesting in the case that its language is
closed under reversal. Infinite words $u_\beta$ associated to
Parry numbers $\beta$ are uniformly recurrent. If $\beta$ is a
simple Parry number, the language of $u_\beta$ is invariant under
reversal if the R\'enyi expansion of 1 satisfies
$d_\beta(1)=tt\cdots ts$, i.e.\ is a confluent Parry number, and
the corresponding palindromic complexity is the subject of this
paper.

For non-simple Parry number $\beta$, the condition under which the
language of the infinite word $u_\beta$ is closed under reversal
has been stated by Bernat~\cite{bernat-thesis}. He has shown that the
language of $u_{\beta}$ is closed under reversal if and only if $\beta$ is a quadratic 
number, i.e.\ a root of minimal polynomial $X^2-aX+b$, with $a \ge b+2$ and $b \ge 1$. In
this case $d_\beta(1)=(a-1) (a-b-1)^\omega$. The palindromic
complexity of the corresponding infinite words $u_\beta$ is
described in~\cite{Lubka}.

Infinite words $u_\beta$ for non-simple Parry numbers $\beta$ are
thus another example for which the equality
$$
{\mathcal P}(n) + {\mathcal P}(n+1) = \Delta \C(n) + 2
$$
is satisfied for all $n\in\N$. According to our knowledge, among
all examples of infinite words satisfying this equality, the words
$u_\beta$ (for both simple and non-simple Parry number $\beta$)
are exceptional in that they have the second difference $\Delta^2
\C(n)\neq0$.

\section*{Acknowledgements}

The authors acknowledge financial support by Czech Science
Foundation GA \v{C}R 201/05/0169, by the grant LC00602 of the
 Ministry of Education, Youth and Sports of the Czech
Republic, and by the program ``Num\'eration" of the ACI {\sl
Nouvelles Interfaces des Math\'ematiques} of the French Ministry
of Research. Part of this work was done while E. P. was Invited
Professor at University Paris~7.


\end{document}